\documentclass[12pt]{amsart}
\usepackage[english]{babel}
\usepackage{amsmath}
\usepackage{amsthm}
\usepackage{amssymb}
\usepackage{mathrsfs}
\usepackage{enumerate}
\usepackage[notcite, final, notref]{showkeys}
\usepackage{dsfont}
\usepackage[dvips]{color}
\usepackage{xcolor}
\newtheorem{theorem}{Theorem}[section]

\newtheorem{lemma}[theorem]{Lemma}

\newtheorem{corol}[theorem]{Corollary}
\newtheorem{definition}[theorem]{Definition}
\theoremstyle{definition}
\newtheorem{remark}[theorem]{Remark}

\def\N{\mathcal{N}}

\def\C{\bf C}
\def\D{\mathcal D}
\def\E{\bf E}
\def\R{\bf R}

\title[]{Higher Dimensional Versions of the Douglas-Ahlfors Identities}

\author[Y. Yang]{Yan Yang}
\address{Yan Yang, School of Mathematics (Zhuhai)\\
  Sun Yat-Sen University\\
 China}
\email{yangyan8@mail.sysu.edu.cn}

\author[T. Qian]{Tao Qian*}
\address{Tao QIAN, Macao Center for Mathematical Sciences\\
Macau University of Science and Technology\\
Macau}
\email{tqian@must.edu.mo}
\thanks{*Corresponding author.\\
Funded by the Science and Technology Development Fund of Macau SAR (grant number 0128/2022/A)}

\begin{document}
\maketitle
\def\e{\bf e}
\begin{abstract}
Denote by ${\mathcal D}$ the open unit disc in the complex plane and $\partial {\mathcal D}$ its boundary. Douglas showed through an identical quantity represented by the Fourier coefficients of the concerned function $u$ that
\begin{eqnarray}\label{abs}
A(u)=\int_{\mathcal D}|\bigtriangledown U|^2dxdy&=&\frac{1}{2\pi}\int\int_{\partial {\mathcal D}\times \partial {\mathcal D}} \left|\frac{u(z_1)-u(z_2)}{z_1-z_2}\right|^2|dz_1||dz_2|,\end{eqnarray}
\end{abstract}
where $u\in L^2(\partial {\mathcal D}), U$ is the harmonic extension of $u$ into
${\mathcal D}$. Ahlfors gave a fourth equivalence form of $A(u)$ in (\ref{more}) via a different proof.
The present article studies relations between the counterpart quantities in higher dimensional spheres with several different but commonly adopted settings, namely, harmonic functions in the Euclidean ${\mathbb R}^n, n\ge 2,$ regular functions in the quaternionic algebra, and Clifford monogenic functions with the real-Clifford algebra ${\mathcal{CL}}_{0, n-1},$ the latter being generated by the multiplication anti-commutative basic imaginary units ${\e}_1, {\e}_2, \cdots , {\e}_{n-1}$ with ${\e}_j^2=-1, j=1, 2, \cdots, n-1.$ It is noted that, while exactly the same equivalence relations hold for harmonic functions in ${\mathbb R}^n$ and regular functions in the quaternionic algebra, for the Clifford algebra setting $n>2,$ the relation (\ref{more}) has to be replaced by essentially a different rule.
\\

\\

\date{today}

\section{Introduction}

In \cite{Douglas} Douglas established the result (\ref{abs}) that was used as the main technical tool to solve the minimum surface problem of Plateau. He proved the identical relation between the two quantities through a third equivalent quantity in terms of the Fourier coefficients of the boundary function $u$, namely $4\sum_{k=1}^\infty k|c_k|^2,$ where $c_k=\frac{1}{2\sqrt{\pi}}\int_0^{2\pi}u(e^{it})e^{-ikt}dt.$ This third quantity is, in fact, closely related to the norm of the Dirichlet space function induced by the  boundary function $u.$

In \cite{Ahlfors}, Ahlfors presented an alternative proof based on differentiation of the Schwarz kernel. Although Ahlfors' proof is not be as direct as Douglas's, it provides a fourth equivalent quantity, namely,
 \begin{eqnarray}\label{more} A(u)=\frac{1}{2}\int_{|z|=1}\overline{f}f'\frac{dz}{i}, \end{eqnarray}
 where $f$ is the holomorphic function obtained through adding the normalized conjugate harmonic part $V$:
\[f(z)= (U+iV)(z)=\frac{1}{2\pi}\int_0^{2\pi}\frac{e^{i\phi}+z}{e^{i\phi}-z}u(\phi)d\phi, \quad V(0)=0.\]
There then followed, besides the application to minimum surface problem, other relevant studies as well, including capacity of sets in planar domains (\cite{Ahlfors}). Enhancing the quantitative equivalence relations in domains other than the unit disc, there recently arouse studies of Wei and Zinsmeister (\cite{WZ}) on domains defined through chord-arc curves and on the related Dirichlet functional spaces. As a very interesting result they show that those quantitative equivalence results stand as a characteristic property of rectifiable Jordan chord-arc curves. The wide and deep connections and applications of the equivalence relations to minimum surfaces, operator and function space theories, as well as in Teichm\"uller theory (see e.g.\cite{NS}) naturally motivate great interest in the related research. This includes the question that in higher dimensions in what extent similar quantities can be defined and whether they have the sort of quantitative equivalence. To our knowledge implementation of the Douglas' program has not been carried out to contexts other than one complex variable. The present paper, as a start, establishes fundamental results in the several real variable cases.  We show that in the unit ball context all the four quantities are well defined and equivalent. Furthermore, we show that there exist counterpart results in the quaternionic and Clifford algebra contexts as well. With the Clifford algebra case there exists an exception (see \S 4).

In the sequel the writing of the paper is divided into three sections. \S 2 deals with the forms of $A(u)$ in the harmonic function setting. Pursuing Douglas's idea, the section proceeds fundamental computations based on the spherical Fourier-Laplace expansions. \S 3 gives the forms of $A(u)$ in the quaternionic setting. The quaternionic setting is different from ${\mathbb R}^4,$ for itself possesses a regular function theory, or a Cauchy-Riemann structure, as in the one-complex variable case. In \S 4 the forms of $A(u)$ in the Clifford algebra setting is presented. The Clifford setting is different from the harmonic, and different from the quaternionic one either, for the domain of the functions in the case is the linear Euclidean space, not an algebra, and the range of the functions, although being an algebra, is not divisible, and nor commutative either. Clifford algebra provides Euclidean space with a Cauchy-Riemann structure as well. The Clifford Cauchy-Riemann equations, however, are not as convenient as for the quaternionic case.  For our purpose some critical computations concerning monogenic functions and the Clifford Schwarz kernel are spelt out that would be new in the foundation of Clifford analysis.

\section {The forms of $A(u)$ in the harmonic analysis setting}

\begin{definition}

For real-valued function $u\in L^2(S^{n-1}), n\ge 2,$ define the functional
\begin{eqnarray}
A(u):= \frac{1}{\omega_{n-1}}\int_{S^{n-1}}\int_{S^{n-1}}\frac{|u(\underline{\eta_1})-u(\underline{\eta_2})|^2}{|\underline{\eta_1}-\underline{\eta_2}|^n}
 dS_{\underline{\eta_1}}dS_{\underline{\eta_2}}.
\end{eqnarray}
\end{definition}

Clearly, the above is an improper double integral that has a determinate positive value, finite or $+\infty$.\\

Let $U$ be the harmonic extension of $u$ into the unit ball $B_{n}$  in ${\mathbb R}^n$. It is well known knowledge that the function values of $u$ are the non-tangential boundary limits of $U$ on the sphere $S^{n-1}$ a.e.  As main result of this section, we have

\begin{theorem}\label{th1}
Let $n\geq 2$. Then

\begin{eqnarray}
A(u)&: =& \frac{1}{\omega_{n-1}}\int_{S^{n-1}}\int_{S^{n-1}}\frac{|u(\underline{\eta_1})-u(\underline{\eta_2})|^2}{|\underline{\eta_1}-\underline{\eta_2}|^n}
 dS_{\underline{\eta_1}}dS_{\underline{\eta_2}}\nonumber\\
& =& \int_{B_n} |\bigtriangledown U|^2 dV\nonumber\\
&=&\sum_{k=1}^{+\infty}k\sum_{j=1}^{a^n_k}|b_j|^2,\nonumber
\end{eqnarray}
where $b_j$'s are the coefficients of the Fourier-Laplace series expansion of $u$ which are given by formula (\ref{eqY14}).

\end{theorem}

When $n=2$, denote by $\D=B_2$ the unit disk, and $T=S^1$ the unit circle. Then there exists the classical result as cited in \cite{WZ}:

\begin{corol}\label{WZ}

\begin{eqnarray}\label{more1}
A(u)&: =&  \frac{1}{2\pi}\int_{T}\int_{T}\frac{|u(z_1)-u(z_2)|^2}{|z_1-z_2|^2}|dz_1| |dz_2|\nonumber\\
&=& \int_{\D} |\bigtriangledown U|^2 dxdy\nonumber\\
&=& \sum_{k=1}^{+\infty}k(a_k^2+b_k^2),
\end{eqnarray}
where $a_k$'s and $b_k$'s are the coefficients of the Fourier series expansion of $u$ given by formula (\ref{eqY9}).
\end{corol}

Based on the Divergence Theorem, we have
\begin{eqnarray}{\label{eq3}}
\int_{B_n} |\bigtriangledown U|^2 dV=\int_{S^{n-1}}U\frac{\partial U}{\partial {\bf n}}dS.
\end{eqnarray}

Next, we perform homogeneous spherical harmonic decomposition on the sphere.

It is well known that $L^2(S^{n-1})=\oplus\sum_{k=0}^{\infty}{\mathcal H}^n_k$,
where ${\mathcal H}^n_k$ is the $a^n_k$-dimensional linear space of all the $k$-spherical harmonics of $n$ variables,
\begin{eqnarray*}
a^n_k=\left\{
\begin{array}{lll}
\vspace{0.2cm}
1, &&\mbox{if } k=0,\\
(n+2k-2)\frac{(n+k-2)!}{k!(n-2)!}, &&\mbox{if } k\in {\N}^+=\{1,2,\cdots\}.
\end{array}
\right.
\end{eqnarray*}

If $u\in L^2(S^{n-1})$, then in the $L^2$ sense, $f$ has the Laplace-Fourier series expansion:
$$u(\underline{\xi})=\sum_{k=0}^{+\infty}Y_k(u)(\underline{\xi}),$$
where $$Y_k(u)(\underline{\xi})=c_{n,k}\int_{S^{n-1}}u(\underline{\eta})P_k^n(\underline{\xi}\cdot\underline{\eta})dS_{\underline{\eta}}$$
stands for the projection of $f$ onto ${\mathcal H}^n_k,$
where $P_k^n$ are the Gegenbauer polynomials and
\begin{eqnarray}\label{eqY4}
c_{n, k}=\frac{1}{\omega_{n-1}}\frac{(n+2k-2)\Gamma(n+k-1)}{(n+k-2)k!\Gamma(n-1)},
\end{eqnarray}\label{eq4}
$\omega_{n-1}=\frac{2\pi^{\frac{n}{2}}}{\Gamma(\frac{n}{2})}$ is the surface area of $S^{n-1}$.\\

For $0<r<1$, there hold the relations
\begin{eqnarray}\label{eq4}
U(r, \underline{\xi})&=&\sum_{k=0}^{+\infty}r^k Y_k(u)(\underline{\xi})\nonumber\\
&=&\sum_{k=0}^{+\infty}r^kc_{n,k}\int_{S^{n-1}}u(\underline{\eta})P_k^n(\underline{\xi}\cdot\underline{\eta})dS_{\underline{\eta}}\\
&=&\int_{S^{n-1}}u(\underline{\eta})
\sum_{k=0}^{+\infty}r^kc_{n,k}P_k^n(\underline{\xi}\cdot\underline{\eta})
dS_{\underline{\eta}}\nonumber.\end{eqnarray}
Since we also have the Poisson kernel integral representation \begin{eqnarray}
U(r,\underline{\xi})=\frac{1}{\omega_{n-1}}\int_{S^{n-1}}u(\underline{\eta})
\frac{1-r^2}{|\underline{\eta}-r\underline{\xi}|^n}dS_{\underline{\eta}},\nonumber
\end{eqnarray}
we have the spherical representation of the Poisson kernel of $B_n$ at $r\underline{\xi}:$
$$P_r(\underline{\eta}, \underline{\xi})=\frac{1}{\omega_{n-1}}\frac{1-r^2}{|\underline{\eta}-r\underline{\xi}|^n}=
\sum_{k=0}^{+\infty}r^kc_{n,k}P_k^n(\underline{\xi}\cdot\underline{\eta}),$$
which, as harmonic by itself, gives rise to the harmonic extension $U$ of $u$ into the unit ball $B_n.$ \\

Before proving Theorem 1, we prove several technical Lemmas.

\begin{lemma}$^{\mbox{\scriptsize \cite{WL}}}$\label{lemma5}
Let $n\geq 2$, $k\in {\mathcal N}^+$. For any orthonormal base $\{y_1, y_2, \cdots, y_{a^n_k}\}$ of ${\mathcal H}^n_k$ and any $\underline{\xi}, \underline{\eta}\in S^{n-1}$, we have
\begin{eqnarray*}
P^n_k(\underline{\xi}\cdot\underline{\eta})=\frac{1}{c_{n, k}}\sum_{j=1}^{a^n_k}\overline{y_j(\underline{\xi})}y_j(\underline{\eta})
=\frac{1}{c_{n, k}}\sum_{j=1}^{a^n_k}\overline{y_j(\underline{\eta})}y_j(\underline{\xi}).
\end{eqnarray*}
\end{lemma}

The proof is simple. The same type expression holds, in fact, for all Hilbert spaces with an orthonormal basis.

Using Lemma \ref{lemma5}, we give another form of the Laplace-Fourier series:

\begin{eqnarray}
u(\underline{\xi})&=&\sum_{k=0}^{+\infty}c_{n,k}\int_{S^{n-1}}u(\underline{\eta})P_k^n(\underline{\xi}\cdot\underline{\eta})dS_{\underline{\eta}}\nonumber\\
&=&\sum_{k=0}^{+\infty}\int_{S^{n-1}}u(\underline{\eta})\sum_{j=1}^{a^n_k}\overline{y_j(\underline{\eta})}y_j(\underline{\xi})dS_{\underline{\eta}}\nonumber\\
&=&\sum_{k=0}^{+\infty}\sum_{j=1}^{a^n_k} b_j y_j(\underline{\xi})\nonumber,
\end{eqnarray}
where
\begin{equation}\label{eqY14}
b_j=\int_{S^{n-1}}u(\underline{\eta})\overline{y_j(\underline{\eta})}dS_{\underline{\eta}}.
\end{equation}

\begin{remark}
Particularly, under the inner product defined by (\ref{eqY14}), when $n=2$, let $\underline{\xi}=(\cos t, \sin t)=e^{it}$ and $\underline{\eta}=(\cos \tau, \sin \tau)=e^{i\tau}$, the orthonormal base of ${\mathcal H}^2_k$ is $$\{\frac{1}{\sqrt{2\pi}}, \frac{1}{\sqrt{\pi}}\cos k t,   \frac{1}{\sqrt{\pi}}\sin kt, \}_{k\in \N^+ }$$ and
$$P_k^2\left(\cos(t-\tau)\right)=\cos k(t-\tau).$$
Then
\begin{equation}\label{eqY9}
u(e^{it})=\frac{a_0}{2}+\sum_{k=1}^{+\infty}(a_k \cos kt+b_k\sin kt).
\end{equation}
\end{remark}

Lemma \ref{lemma5} immediately implies
\begin{lemma}$^{\mbox{\scriptsize \cite{WL}}}$ (Funk-Hecke formula)\\

If $g\in {\mathcal H}_k^n$, for any $\underline{\xi}\in S^{n-1}$, we have
$$c_{n, k}\int_{S^{n-1}}P_k^n(\underline{\xi}\cdot \underline{\eta})g(\underline{\eta})dS_{\underline{\eta}}=g(\underline{\xi}).$$
\end{lemma}

\begin{lemma}$^{\mbox{\scriptsize \cite{WL}}}$\label{lem3}

Let $-1<t<1.$ When $n>2$,
$$\sum_{k=0}^{+\infty}\frac{\Gamma(n+k-2)}{k!\Gamma(n-2)}x^k P^n_k(t)
=\frac{1}{(1-2tx+x^2)^{\frac{n-2}{2}}}.$$

When $n=2$,
$$\sum_{k=1}^{+\infty}\frac{1}{k}x^k P^2_k(t)
=-\frac{1}{2}\ln(1-2xt+x^2).$$

\end{lemma}

\begin{lemma}\label{lem4}
Let $-1<t<1$, $n\geq 2$, then
\begin{eqnarray}\label{eqY7}
&&\sum_{k=1}^{+\infty}kr^{2k+n-2}c_{n,k}P^n_k(t)\nonumber\\
&=&\frac{1}{\omega_{n-1}}r^n\mbox{ }\frac{nt-(n-4)tr^4-(n+2)r^2+(n-2)r^6}{(1-2tr^2+r^4)^{\frac{n}{2}+1}}.
\end{eqnarray}
We denote $$J(r, t)=-r^n\mbox{ }\frac{nt-(n-4)tr^4-(n+2)r^2+(n-2)r^6}{(1-2tr^2+r^4)^{\frac{n}{2}+1}}.$$
\end{lemma}

{\bf Proof: }
When $n>2,$ taking derivative with respect to $x$ on both sides of the equality in Lemma \ref{lem3}, we have
$$\sum_{k=1}^{+\infty}\frac{\Gamma(n+k-2)}{k!\Gamma(n-2)}k x^{k-1} P^n_k(t)
=\frac{(n-2)(t-x)}{(1-2tx+x^2)^{\frac{n}{2}}}.$$

Then $$\sum_{k=1}^{+\infty}\frac{\Gamma(n+k-2)}{(k-1)!\Gamma(n-1)} x^{k-1} P^n_k(t)
=\frac{t-x}{(1-2tx+x^2)^{\frac{n}{2}}}.$$

Letting $x=r^2$ and multiplying $r^n$ to both sides, we obtain

$$
\sum_{k=1}^{+\infty}\frac{\Gamma(n+k-2)}{(k-1)!\Gamma(n-1)} r^{n+2k-2} P^n_k(t)
=\frac{tr^n-r^{n+2}}{(1-2tr^2+r^2)^{\frac{n}{2}}}.$$
Taking derivative with respect to $r$ on both sides, we have

$$
\sum_{k=1}^{+\infty}\frac{\Gamma(n+k-2)}{(k-1)!\Gamma(n-1)}(n+2k-2) r^{n+2k-3} P^n_k(t)
=\frac{d}{dr}\left(\frac{tr^n-r^{n+2}}{(1-2tr^2+r^2)^{\frac{n}{2}}}\right).$$

By replacing the expression of $c_{n, k}$ in equation (\ref{eqY4}), we get
\begin{eqnarray}
&&\sum_{k=1}^{+\infty}kr^{n+2k-2}c_{n,k}P^n_k(t)\nonumber\\
&=&\frac{1}{\omega_{n-1}}\sum_{k=1}^{+\infty}\frac{\Gamma(n+k-1)(n+2k-2)}{(n+k-2)(k-1)!\Gamma(n-1)} r^{n+2k-2} P^n_k(t)\nonumber\\
&=&\frac{1}{\omega_{n-1}}r\frac{d}{dr}\left(\frac{tr^n-r^{n+2}}{(1-2tr^2+r^2)^{\frac{n}{2}}}\right)\nonumber\\
&=&\frac{1}{ \omega_{n-1}}r^n\mbox{ }\frac{nt-(n-4)tr^4-(n+2)r^2+(n-2)r^6}{(1-2tr^2+r^4)^{\frac{n}{2}+1}}\nonumber.
\end{eqnarray}

When $n=2$.
$$\sum_{k=1}^{+\infty}\frac{1}{k}x^k P^2_k(t)
=-\frac{1}{2}\ln(1-2xt+x^2).$$

Using a similar method, we have

$$\sum_{k=1}^{+\infty}kx^{2k} P^2_k(t)
=r^2\frac{t+tr^4-2r^2}{(1-2tr^2+r^4)^2}.$$

In short, when $n\geq2$, we have formula (\ref{eqY7}). The proof is complete.

\bigskip
{\bf Proof of Theorem \ref{th1}: } By Levi's Theorem and formula $(\ref{eq3})$, we have
\begin{eqnarray}\label{eq6}
\int_{B_n} |\bigtriangledown U|^2 dV&=&\lim_{r\rightarrow 1^{-}}\int_{B_n(r)} |\bigtriangledown U|^2 dV\nonumber\\
&=&\lim_{r \rightarrow 1^{-}}\int_{S_r^{n-1}}U\frac{\partial U}{\partial {\bf n}}dS\nonumber\\
&=&\lim_{r \rightarrow 1^{-}}\int_{S^{n-1}}U\frac{\partial U}{\partial {r}}r^{n-1}dS.
\end{eqnarray}
Using formula ($\ref{eq4}$), we obtain
$$U=c_{n, 0}\int_{S^{n-1}}u(\underline{\eta})dS_{\underline{\eta}}+\sum_{k=1}^{+\infty}r^kc_{n,k}\int_{S^{n-1}}u(\underline{\eta})P_k^n(\underline{\xi}\cdot\underline{\eta})dS_{\underline{\eta}}.$$
Then $$\frac{\partial U}{\partial {r}}r^{n-1}=\sum_{k=1}^{+\infty}kr^{k+n-2}c_{n,k}\int_{S^{n-1}}u(\underline{\eta})P_k^n(\underline{\xi}\cdot\underline{\eta})dS_{\underline{\eta}}.$$
Therefore, we get
\begin{eqnarray*}
\int_{S^{n-1}}U\frac{\partial U}{\partial {r}}r^{n-1}dS
&=&\sum_{k=1}^{+\infty}kr^{2k+n-2}c^2_{n,k}\int_{S^{n-1}}\int_{S^{n-1}}u(\underline{\eta_1})u(\underline{\eta_2})\\
&&\left[\int_{S^{n-1}}P^n_k (\underline{\xi}\cdot\underline{\eta_1})P^n_k(\underline{\xi}\cdot\underline{\eta_2})dS_{\underline{\xi}}\right]dS_{\underline{\eta_1}}dS_{\underline{\eta_2}},
\end{eqnarray*}
where we used orthogonality between $P^n_k(\underline{\xi}\cdot\underline{\eta})$ and  $P^n_j(\underline{\xi}\cdot\underline{\eta})$ when $k\neq j$.

By the Funk-Hecke formula, we have
$$c_{n, k}\int_{S^{n-1}}
P^n_k (\underline{\xi}\cdot\underline{\eta_1})P^n_k(\underline{\xi}\cdot\underline{\eta_2})dS_{\underline{\xi}}= P^n_k(\underline{\eta_1}\cdot\underline{\eta_2}).$$

Then, we obtain

\begin{eqnarray}\label{Yeq10}
&&\int_{S^{n-1}}U\frac{\partial U}{\partial {r}}r^{n-1}dS\nonumber\\
&=&\sum_{k=1}^{+\infty}kr^{2k+n-2}c_{n,k}\int_{S^{n-1}}\int_{S^{n-1}}u(\underline{\eta_1})u(\underline{\eta_2})P^n_k(\underline{\eta_1}\cdot\underline{\eta_2})dS_{\underline{\eta_1}}dS_{\underline{\eta_2}}.
\end{eqnarray}

On the one hand, from formula (\ref{Yeq10}), Lemma \ref{lemma5} and formula (\ref{eqY14}), we have

\begin{eqnarray}\label{eqyy15}
&&\int_{S^{n-1}}U\frac{\partial U}{\partial {r}}r^{n-1}dS\nonumber\\
&=&\sum_{k=1}^{+\infty}kr^{2k+n-2}c_{n,k}\int_{S^{n-1}}\int_{S^{n-1}}u(\underline{\eta_1})u(\underline{\eta_2})P^n_k(\underline{\eta_1}\cdot\underline{\eta_2})dS_{\underline{\eta_1}}dS_{\underline{\eta_2}}\nonumber\\
&=&\sum_{k=1}^{+\infty}kr^{2k+n-2}\int_{S^{n-1}}\int_{S^{n-1}}u(\underline{\eta_1})u(\underline{\eta_2})\sum_{j=1}^{a^n_k}\overline{y_j(\underline{\eta_1})}y_j(\underline{\eta_2})dS_{\underline{\eta_1}}dS_{\underline{\eta_2}}\nonumber\\
&=&\sum_{k=1}^{+\infty}kr^{2k+n-2}\sum_{j=1}^{a^n_k}\int_{S^{n-1}}u(\underline{\eta_1})\overline{y_j(\underline{\eta_1})}dS_{\underline{\eta_1}}\int_{S^{n-1}}\overline{u(\underline{\eta_2})\overline{y_j(\underline{\eta_2})}}dS_{\underline{\eta_2}}\nonumber\\
&=&\sum_{k=1}^{+\infty}kr^{2k+n-2}\sum_{j=1}^{a^n_k}|b_j|^2.
\end{eqnarray}

Therefore, by formula (\ref{eq6}) and (\ref{eqyy15}), we obtain
\begin{eqnarray}\label{eq16}
\int_{B_n} |\bigtriangledown U|^2 dV&=&\lim_{r\rightarrow 1^{-}}\int_{B_n(r)} |\bigtriangledown U|^2 dV\nonumber\\
&=&\lim_{r \rightarrow 1^{-}}\int_{S^{n-1}}U\frac{\partial U}{\partial {r}}r^{n-1}dS\nonumber\\
&=&\lim_{r \rightarrow 1^{-}}\sum_{k=1}^{+\infty}kr^{2k+n-2}\sum_{j=1}^{a^n_k}|b_j|^2\nonumber.
\end{eqnarray}
If $$\sum_{k=1}^{+\infty}k\sum_{j=1}^{a^n_k}|b_j|^2<+\infty,$$  using the second theorem of Abel, we have
$$\lim_{r \rightarrow 1^{-}}\sum_{k=1}^{+\infty}kr^{2k+n-2}\sum_{j=1}^{a^n_k}|b_j|^2=\sum_{k=1}^{+\infty}k\sum_{j=1}^{a^n_k}|b_j|^2.$$
Then
\begin{eqnarray}\label{eq15}
&&\int_{S^{n-1}}U\frac{\partial U}{\partial {r}}r^{n-1}dS\nonumber\\
&=&\lim_{r \rightarrow 1^{-}}\sum_{k=1}^{+\infty}kr^{2k+n-2}\sum_{j=1}^{a^n_k}|b_j|^2\nonumber\\
&=&\sum_{k=1}^{+\infty}k\sum_{j=1}^{a^n_k}|b_j|^2.\nonumber
\end{eqnarray}

If $$\sum_{k=1}^{+\infty}k\sum_{j=1}^{a^n_k}|b_j|^2=+\infty,$$  using the inverse theorem of the second theorem of Abel, we have
$$\lim_{r \rightarrow 1^{-}}\sum_{k=1}^{+\infty}kr^{2k+n-2}\sum_{j=1}^{a^n_k}|b_j|^2=+\infty.$$

Therefore, we get
$$\int_{B_n} |\bigtriangledown U|^2 dV=\lim_{r \rightarrow 1^{-}}\sum_{k=1}^{+\infty}kr^{2k+n-2}\sum_{j=1}^{a^n_k}|b_j|^2=+\infty.$$

On the other hand,
\begin{eqnarray}
&&\int_{S^{n-1}}U\frac{\partial U}{\partial {r}}r^{n-1}dS\nonumber\\
&=&-\frac{1}{2}\sum_{k=1}^{+\infty}kr^{2k+n-2}c_{n,k}\int_{S^{n-1}}\int_{S^{n-1}}[u(\underline{\eta_1})-u(\underline{\eta_2})]^2P^n_k(\underline{\eta_1}\cdot\underline{\eta_2})dS_{\underline{\eta_1}}dS_{\underline{\eta_2}}\nonumber\\
&=&\int_{S^{n-1}}\int_{S^{n-1}}[u(\underline{\eta_1})-u(\underline{\eta_2})]^2[-\frac{1}{2}\sum_{k=1}^{+\infty}kr^{2k+n-2}c_{n,k}P^n_k(\underline{\eta_1}\cdot\underline{\eta_2})]dS_{\underline{\eta_1}}dS_{\underline{\eta_2}}\nonumber.
\end{eqnarray}

Using Lemma \ref{lem4}, we have
\begin{eqnarray}\label{eq7}
&&\int_{S^{n-1}}U\frac{\partial U}{\partial {r}}r^{n-1}dS\nonumber\\
&=&\frac{1}{2 \omega_{n-1}}\int_{S^{n-1}}\int_{S^{n-1}}[u(\underline{\eta_1})-u(\underline{\eta_2})]^2 J(r,\cos\theta)dS_{\underline{\eta_1}}dS_{\underline{\eta_2}},
\end{eqnarray}
where $\theta$ is the angle between $\underline{\eta_1}$ and $\underline{\eta_2}$, then $\underline{\eta_1}\cdot\underline{\eta_2}=\cos\theta$ and $|\underline{\eta_1}-\underline{\eta_2}|^2=4\sin^2\frac{\theta}{2}$.
By calculating directly, we have
\begin{eqnarray*}
&&J(r, \cos\theta)\\
&=& r^n \frac{n\left[(1+r^2)^2\sin^2\frac{\theta}{2}-(1-r^2)^2\cos^2\frac{\theta}{2}\right]-(n-2)r^2\left[(1+r^2)^2\sin^2\frac{\theta}{2}+(1-r^2)^2\cos^2\frac{\theta}{2}\right]}{[(1+r^2)^2\sin^2\frac{\theta}{2}+(1-r^2)^2\cos^2\frac{\theta}{2}]^{\frac{n}{2}+1}}\\
&\rightarrow& \frac{2}{2^n (\sin^2\frac{\theta}{2})^{\frac{n}{2}}}=J(1, \cos\theta) \mbox{  }\mbox{ }(r\rightarrow 1^-).
\end{eqnarray*}
In particular, when $r$ is close to 1, $J(r, \cos\theta)$ is positive.

If
\begin{eqnarray}\label{eq8}
&&\lim_{r\rightarrow 1^-}\int_{S^{n-1}}\int_{S^{n-1}}[u(\underline{\eta_1})-u(\underline{\eta_2})]^2J(r, \cos\theta)dS_{\underline{\eta_1}}dS_{\underline{\eta_2}}\nonumber\\
&=&\int_{S^{n-1}}\int_{S^{n-1}}[u(\underline{\eta_1})-u(\underline{\eta_2})]^2 J(1, \cos\theta)dS_{\underline{\eta_1}}dS_{\underline{\eta_2}},
\end{eqnarray}
then by formula (\ref{eq6}), (\ref{eq7}) and (\ref{eq8}), we get
\begin{eqnarray*}
&&\int_{B_n} |\bigtriangledown U|^2 dV\\
&=&\lim_{r\rightarrow 1^-}\frac{1}{2\omega_{n-1}}\int_{S^{n-1}}\int_{S^{n-1}}[u(\underline{\eta_1})-u(\underline{\eta_2})]^2 J(r, \cos\theta)dS_{\underline{\eta_1}}dS_{\underline{\eta_2}}\\
&=&\frac{1}{2\omega_{n-1}}\int_{S^{n-1}}\int_{S^{n-1}}[u(\underline{\eta_1})-u(\underline{\eta_2})]^2J(1, \cos\theta)dS_{\underline{\eta_1}}dS_{\underline{\eta_2}}\\
&=&\frac{1}{\omega_{n-1}}\int_{S^{n-1}}\int_{S^{n-1}}\frac{|u(\underline{\eta_1})-u(\underline{\eta_2})|^2}{|\underline{\eta_1}-\underline{\eta_2}|^n}
 dS_{\underline{\eta_1}}dS_{\underline{\eta_2}}.
\end{eqnarray*}

In fact,
\begin{eqnarray*}
&&\frac{|J(r, \cos\theta)|}{J(1, \cos\theta)}\\
&=&\frac{1}{2} \frac{2^n (\sin^2\frac{\theta}{2})^{\frac{n}{2}}r^n}{[(1+r^2)^2\sin^2\frac{\theta}{2}+(1-r^2)^2\cos^2\frac{\theta}{2}]^{\frac{n}{2}}}\\
&\times&\frac{\left|-n\left[(1+r^2)^2\sin^2\frac{\theta}{2}-(1-r^2)^2\cos^2\frac{\theta}{2}\right]+(n-2)r^2\left[(1+r^2)^2\sin^2\frac{\theta}{2}+(1-r^2)^2\cos^2\frac{\theta}{2}\right]\right|}{(1+r^2)^2\sin^2\frac{\theta}{2}+(1-r^2)^2\cos^2\frac{\theta}{2}}\\
&\leq&\frac{1}{2}[n+(n-2)r^2]\\
&\leq&\frac{1}{2}[n+(n-2)]=n-1.
\end{eqnarray*}

If $$\int_{S^{n-1}}\int_{S^{n-1}}[u(\underline{\eta_1})-u(\underline{\eta_2})]^2 J(1, \cos\theta)dS_{\underline{\eta_1}}dS_{\underline{\eta_2}}< +\infty,$$ then using the Lebesgue dominated convergence theorem, formula (\ref{eq8}) holds.

If $$\int_{S^{n-1}}\int_{S^{n-1}}[u(\underline{\eta_1})-u(\underline{\eta_2})]^2 J(1, \cos\theta)dS_{\underline{\eta_1}}dS_{\underline{\eta_2}}= +\infty,$$
Using Fatou's Lemma, we have
\begin{eqnarray*}
&&\lim_{r\to 1^-}\int_{S^{n-1}}\int_{S^{n-1}}[u(\underline{\eta_1})-u(\underline{\eta_2})]^2J(r, \cos\theta)dS_{\underline{\eta_1}}dS_{\underline{\eta_2}}\\
&\geq& \int_{S^{n-1}}\int_{S^{n-1}}[u(\underline{\eta_1})-u(\underline{\eta_2})]^2 J(1, \cos\theta)dS_{\underline{\eta_1}}dS_{\underline{\eta_2}}\\
&=&+\infty.
\end{eqnarray*}
Therefore, formula (\ref{eq8}) also holds. The proof of Theorem \ref{th1} is complete.

\bigskip

\section{The forms of $A(u)$ in the Quaternionic analysis setting}

\def\H{\mathbb{H}}
\def\HC{{H}^2(\mathbb{C}^+)}
\def\R{\mathbb{R}}
\def\N{\mathbb{N}}
\def\C{\mathbb{C}}
\def\e{\bf e}
\def\E{\bf E}
\def\i{\bf i}
\def\j{\bf j}
\def\k{\bf k}
It is well known that if $U$ is harmonic in $\D$, then there exists a canonical conjugate harmonic function $V$ in $\D$, such that $F=U+iV$ is holomorphic in $\D$.

Using the Cauchy-Riemann equation, we have
\begin{eqnarray*}
\left\{
\begin{array}{lll}
\vspace{0.2cm}
\frac{\partial U}{\partial x}=\frac{\partial V}{\partial y}\\

\frac{\partial U}{\partial y}=-\frac{\partial V}{\partial x}.
\end{array}
\right.
\end{eqnarray*}

Therefore,
$$F'(z)=\frac{\partial U}{\partial x}+i\frac{\partial V}{\partial x}= \frac{\partial U}{\partial x}-i\frac{\partial U}{\partial y}.$$
Then there exists another form of $A(f)$ in 2-dimensional case, namely

$$A(u)=\int_{\D} |\bigtriangledown U|^2 dxdy=\int_{\D} |F'(z)|^2 dxdy,$$
which is essentially (\ref{more}) through using Green's formula.\\

In the quaternionic four dimensional case, we will be able to get a similar counterpart form, that is, specifically,

\begin{eqnarray}\label{eqyy21}
A(u) = \int_{B_{4}} |\bigtriangledown U|^2 dV=\int_{B_4} |\bar{D}F|^2 dV,
\end{eqnarray}
where $F$ is the associated quaternionic regular function. Also see the reasoning given after Lemma \ref{lem13}.\\

To deduce the above formula we recall some basic knowledge about quaternion algebra and quaternionic analysis. For more details, see \cite{Sud}.

The real quaternions constitute a four-dimensional algebra which is non-commutative but associative and divisible.\\

We choose a basis $1, \mbox{ }\i, \mbox{ }\j, \mbox{ }\k$ such that the multiplication is given by the rules as follows:
$$\i^2=\j^2=\k^2=-1,$$
$$\i \j=-\j \i=\k, \mbox{ }\j \k=-\k \j=\i,\mbox{ } \k \i=-\i \k=\j.$$

A typical quaternion is denoted as $$q: = x_0+x_1{\i}+x_2{\j}+x_3{\k}, $$ where $x_k\in {\R}, k=0, 1, 2, 3.$

Denote $${\H}: =\{q\ |\ q=x_0+x_1{\i}+x_2{\j}+x_3{\k}, x_k\in {\R}, k=0,1,2,3\}$$ the set of quaternions.

Let $q=x_0+x_1{\i}+x_2{\j}+x_3{\k}\in \H$, define $\bar{q}: = x_0-x_1{\i}-x_2{\j}-x_3{\k}$ to be the conjugate of $q$. Then we have
$$q\bar{q}=\bar{q}q=\sum_{k=0}^{3}x_k^2.$$ We define $|q|: =\sqrt{\sum_{k=0}^{3}x_k^2}$.

\begin{definition}
For any $q\in \H$, $q= x_0+x_1{\i}+x_2{\j}+x_3{\k}$, we call $x_0$ the real part of $q$, also denoted as $Re(q),$ and $x_1{\i}+x_2{\j}+x_3{\k}$ the non-real part of $q$, denoted as ${\rm {NRe}}(q)$.
\end{definition}

Define the Cauchy-Riemann-Fueter operator
$$D=\frac{1}{2}(\frac{\partial}{\partial x_0}+\frac{\partial}{\partial x_1}{\i}+\frac{\partial}{\partial x_2}{\j}+\frac{\partial}{\partial x_3}{\k})$$
and its conjugate operator
$$\bar{D}=\frac{1}{2}(\frac{\partial}{\partial x_0}-\frac{\partial}{\partial x_1}{\i}-\frac{\partial}{\partial x_2}{\j}-\frac{\partial}{\partial x_3}{\k}).$$
Clearly, the Laplace operator in the 4-D space $\Delta_q=4D\bar{D}=4\bar{D}D$.

Let $g(q): {\H}\rightarrow {\H}$ be a quaternionic valued function.

Denote $g(q)=g_0(q)+g_1(q){\i}+g_2(q){\j}+g_3(q){\k}$, where $g_k(q), k=0, 1, 2, 3$ are  real-valued functions.

\begin{definition}
We call $g(q)$ left-regular in domain $\Omega\subseteq{\H}$, if $g\in C^1(\Omega)$ and $Dg=0$. If $g$ is left-regular, then we call $\bar{D}g$ the left-derivative of $g$.
\end{definition}

\begin{remark}

(1). If $g$ is left-regular, then

\begin{equation}\label{eqyy22}
\bar{D}g=\frac{\partial g}{\partial x_0}.
\end{equation}

(2). If $g$ is left-regular in $\Omega$,  then $g,$ and each $g_i, i=0,1,2,3,$ as well, is harmonic in $\Omega$.
\end{remark}

\begin{lemma}$^{\mbox{\scriptsize \cite{Sud}}}$

\begin{eqnarray}\label{eq24}
Dg=0 \Longleftrightarrow \left\{
\begin{array}{lll}
\vspace{0.2cm}
\displaystyle\frac{\partial g_0}{\partial x_0}=\sum_{k=1}^3\frac{\partial g_k}{\partial x_k} \\
\vspace{0.2cm}
\displaystyle\frac{\partial g_0}{\partial x_k}=-\frac{\partial g_k}{\partial x_0}, k=1, 2, 3 \\
\vspace{0.2cm}
\displaystyle\frac{\partial g_k}{\partial x_j}=-\frac{\partial g_j}{\partial x_k}, j \not= k.
\end{array}
\right.
\end{eqnarray}

\end{lemma}

\begin{lemma}$^{\mbox{\scriptsize \cite{Sud}}}$\label{lem13}

Let $U$ be a real-valued function defined on a star-shaped open set $\Omega\subseteq {\H}$. If $U$ is harmonic, then the unique left-regular function $F$ on $\Omega$ such that $Re F=U$ and normalized by $F(0)=U(0)$ is given by
$$F(q)=U(q)+2 {\rm{NRe}}\int_0^1 s^2\bar{D}U(sq)q ds.$$
\end{lemma}
\bigskip

In general, denote by $F(q)=U(q)+U_1(q){\i}+U_2(q) {\j}+U_3(q) {\k}$ a left-regular function on the ball whose real part is the harmonic function $U$ and the non-real part is normalized to be zero at the origin. Using the formulas in (\ref{eq24}), we have $\frac{\partial U}{\partial x_k}=-\frac{\partial U_k}{\partial x_0}, k=1, 2, 3.$

By virtue of  (\ref{eqyy22}) we have
$$|\bigtriangledown U|^2=\sum_{k=0}^{3}(\frac{\partial U}{\partial x_k})^2=(\frac{\partial U}{\partial x_0})^2+\sum_{k=1}^{3}(\frac{\partial U_k}{\partial x_0})^2=|\bar{D}F|^2.$$

Therefore, (\ref{eqyy21}) holds.

Furthermore, in the quaternion algebra setting, we also can obtain a form of $A(u)$ similar to (\ref{more}) given by Alfors for one complex variable.
That is
\begin{equation}\label{addeq1}
A(u)=\frac{1}{2}\int_{S^3}\bar{F}d\sigma {\bar D}F,
\end{equation}
where
$$d\sigma=dx_1\wedge dx_2\wedge dx_3-{\i} dx_0\wedge dx_2\wedge dx_3+{\j} dx_0\wedge dx_1\wedge dx_3 -{\k}dx_0\wedge dx_1\wedge dx_2.$$

\begin{lemma}(Stokes' formula)

If $g\in C^1(\bar{\Omega})$, then $$\int_{\Omega}dw=\int_{\partial \Omega}w.$$
\end{lemma}

\begin{theorem}
Let $F(x)\in C^1(\bar{\Omega})$. If $F(q)$ is left-regular in $\Omega\subseteq \H$, then we have
\begin{equation}\label{Yy24}
\int_{\Omega}|\bar{D}F|^2 dV=\frac{1}{2}\int_{\partial \Omega}\bar{F}d\sigma{\bar D}F,
\end{equation}
where $dV=dx_0\wedge dx_1\wedge dx_2\wedge dx_3$ is the volume element in $\H$.
\end{theorem}

{\bf Proof: } Using Stokes' formula, we have
$$\int_{\partial \Omega}\bar{F}d\sigma{\bar D}F=\int_{\Omega} d(\bar{F}d\sigma{\bar D}F).$$
Then,
\begin{eqnarray*}
d(\bar{F}d\sigma{\bar D}F)&=&d(\bar{F}d\sigma){\bar D}F+\bar{F}d(d\sigma{\bar D}F)\\
&=&d\bar{F}\wedge d\sigma {\bar D}F-\bar{F}d\sigma\wedge d({\bar D}F)\\
&=&\left[2(\bar{F}D){\bar D}F+2\bar{F}(D{\bar D}F)\right]dV\\
&=&2\left[(\bar{F}D){\bar D}F+\bar{F}(\frac{1}{4}\Delta_q F)\right]dV\\
&=&2|{\bar D}F|^2dV.
\end{eqnarray*}

This completes the proof.\\

Using formula (\ref{eqyy21}) and (\ref{Yy24}), we obtain (\ref{addeq1}).\\

In summary, with the quaternionic analysis setting, we have:

\begin{theorem}\label{th11}
Let $F$ be left-regular in $B_4$ and $U$ its real part with square-integrable boundary value $u$ in $S^3.$ Then there holds
\begin{eqnarray*}
A(u)&:=&\frac{1}{2\pi^2}\int_{S^3}\int_{S^3}\frac{|u(q_1)-u(q_2)|^2}{|q_1-q_2|^4}
 dS_{q_1}dS_{q_2}\\
&=&\sum_{k=1}^{+\infty}k\sum_{j=1}^{a^4_k}|b_j|^2\\
&=&\int_{B_4}|\bigtriangledown U|dV\\
&=&\int_{B_4}|\bar{D}F|^2 dV\\
&=&\frac{1}{2}\int_{S^3}\bar{F}d\sigma \bar{D}F,
\end{eqnarray*}
where $b_j$'s are the coefficients of the Fourier-Laplace series of $u$ which are given by formula (\ref{eqY14}) for $n=4$.
\end{theorem}

\bigskip

\section{Forms of $A(u)$ in the Clifford analysis setting}

The quaternionic four-dimensional case is different from the case where functions are defined on the span of the basic elements ${\e}_0=1, {\e}_1$ and $ {\e}_2$ and take values in the whole algebra generated by those basic elements. We regard the latter as the standard Clifford analysis setting in which functions are defined in the Euclidean space but with one more dimension the range of the functions possesses a Cauchy-Riemann structure. In the case we will show that, unlike the complex and the quaternionic case, in the Clifford analysis setting, when $n>2$ we have
\begin{eqnarray}\label{eq23}
A(u) = \int_{B_{n}} |\bigtriangledown U|^2 dV\not=\int_{B_n} |\bar{D}F|^2 dV=\frac{1}{2}{\rm Sc}\int_{S^{n-1}}\bar{F}d\sigma{\bar D}F.
\end{eqnarray}

Clifford analysis has close connections with harmonic analysis. For instance, a harmonic conjugate system (\cite{SW}) is just the components of a Clifford monogenic function, and Hilbert transform of a function in Euclidean space is the ${\bf e}_j$-multiple-sum of the Riesz transforms (\cite{DMQ}). Clifford algebras structure has been well adopted into contemporary harmonic analysis (\cite{GM}, \cite{LMcQ}, \cite{LMcS}), and helps to solve deep analysis problems.\\

Next, we will introduce some basic knowledge about Clifford algebra and Clifford analysis. For details, see \cite{DSS}.

Let ${\e}_1, {\e}_2, \cdots , {\e}_{n-1}$ be the basic elements satisfying ${\e}_j{\e}_k+{\e}_k{\e}_j=-2\delta_{jk}$, where $\delta_{jk}=1$ if $j=k$, and $\delta_{jk}=0$ otherwise, $j, k=1, 2, \cdots , n-1$. Let
$${\R}^n=\{x\ |\ x=x_0+x_1{\e}_1+\cdots +x_{n-1}{\e}_{n-1}:x_k \in {\R}, k=1, 2, \cdots , {n-1}\}$$
be the n-dimensional non-homogeneous Euclidean space.

Let ${\mathcal{CL}}_{0, n-1}$ denote the real
Clifford algebra generated by ${\e}_1, {\e}_2, \cdots , {\e}_{n-1}$. The linear basis for the
Clifford algebra is given by ${\e}_A$, where $A$ runs over all the
ordered subsets of $\{0, 1, \cdots, n-1\}$, namely,
 $$A=\{1\leq i_1<i_2<\cdots <i_l\leq n-1\}, \mbox{ }1\leq l\leq n-1. $$
We identify ${\e}_0={\e}_{\emptyset}=1.$

A general element $x$ of ${\mathcal{CL}}_{0, n-1}$ can be represented in the form $x=\sum_{k=0}^{n-1}[x]_k$, where
$[x]_k=\sum\limits_{A} x_A {\e}_A, {\e}_A={\e}_{i_1}{\e}_{i_2}\cdots {\e}_{i_k}, 1\leq i_1< i_2<\cdots< i_k\leq {n-1}.$

When $x=\sum_{k=0}^{n-1}[x]_k \in {\mathcal{CL}}_{0, n-1}$, then x consists of a scalar part and a non-scalar part,
denoted, respectively, by
$$x_0={\rm {Sc}}(x), \quad \sum_{k=1}^{n-1}[x]_k={\rm{NSc}}(x). $$

We defien the norm of $x\in {\mathcal{CL}}_{0, n-1}$ to be $|x|=\sqrt{{\rm Sc}(x\bar{x})}=\sqrt{{\rm Sc}(\bar{x}x)}=\sqrt{\sum_{k=0}^{n-1}\sum_A x_A^2}$.

Define the generalized Cauchy-Riemann operator as
$$D=\frac{1}{2}(\frac{\partial}{\partial x_0}+\sum_{k=1}^{n-1}\frac{\partial}{\partial x_k}{\e}_k)$$
and the related conjugate operator
$$\bar{D}=\frac{1}{2}(\frac{\partial}{\partial x_0}-\sum_{k=1}^{n-1}\frac{\partial}{\partial x_k}{\e}_k).$$
Clearly, the Laplace operator $\Delta_x$ in the n-D space satisfies $\Delta_x=4D\bar{D}=4\bar{D}D$.

We will involve the radial decomposition of $D$. Let, as in \cite{DSS}, $x=r\xi$ with $r=|x|$ and $\xi\in S^{n-1}.$ Then $D$ can be expressed as
\begin{equation}
D=\frac{1}{2}(\xi\partial_r+\frac{1}{r}\partial_{\xi})=\frac{1}{2}\xi({\partial_r}+\Gamma_{\xi}),\nonumber
\end{equation}
where $\Gamma_{\xi}=\bar{\xi}\partial_{\xi}$.

\begin{equation}\label{yyeqadd1}
{\bar D}=\frac{1}{2}({\partial_r}+\Gamma^*_{\xi}){\bar \xi},
\end{equation}
where $\Gamma^*_{\xi}=\bar{\partial_{\xi}}\xi$. We note that the action of $\Gamma^*_{\xi}$ on a function $f$ is proceeded as $\Gamma^*_{\xi}f=\bar{\partial_{\xi}}(\xi f)$.\\

As is well known, the Laplacian
\begin{equation}\label{addyyeq2}
\Delta_x=\partial^2_r+\frac{n-1}{r}\partial_r+\frac{1}{r^2}\Delta_{\xi},
\end{equation}
where $\Delta_{\xi}$ is the Laplace-Beltrami operaor on $S^{n-1}$.

While we also have
\begin{eqnarray*}
\Delta_x&=&4\bar{D}D\\
&=&(\partial_r+\frac{1}{r}\Gamma^*_{\xi})\bar{\xi}\xi(\partial_r+\frac{1}{r}\Gamma_{\xi})\\
&=&(\partial_r+\frac{1}{r}\Gamma^*_{\xi})(\partial_r+\frac{1}{r}\Gamma_{\xi})\\
&=&\partial^2_r+\frac{\partial_r}{r}(\Gamma^*_{\xi}+\Gamma_{\xi})+\frac{1}{r^2}(\Gamma^*_{\xi}\Gamma_{\xi}-\Gamma_{\xi}).
\end{eqnarray*}

With formula (\ref{addyyeq2}), we obtain
\begin{equation}\label{addyyeq3}
\Gamma^*_{\xi}+\Gamma_{\xi}=(n-1){\bf I},
\end{equation}
where ${\bf I}$ is the identity operator.\\

\begin{remark}
For the one complex variable case treated in \cite{Douglas} and \cite{Ahlfors} the corresponding objects are $z=x+iy=r\xi$ with $\xi=e^{i\theta}$. In the case we have $$D=\frac{1}{2}(\frac{\partial}{\partial x}+i\frac{\partial}{\partial y})=\frac{1}{2}\xi({\partial_r}+\Gamma_{\xi})
=\frac{1}{2}e^{i\theta}(\frac{\partial}{\partial r}+i\frac{\partial}{\partial \theta})$$ and
$$\bar{D}=\frac{1}{2}(\frac{\partial}{\partial x}-i\frac{\partial}{\partial y})=\frac{1}{2}\bar{\xi}({\partial_r}+\Gamma^*_{\xi})=\frac{1}{2}e^{-i\theta}(\frac{\partial}{\partial r}-i\frac{\partial}{\partial \theta}).$$
We hence have $\Gamma_{\xi}=i\frac{\partial}{\partial \theta}$ and $\Gamma^*_{\xi}=-i\frac{\partial}{\partial \theta},$ and  $\Gamma_{\xi}+\Gamma^*_{\xi}=0$. The above can be interpreted as the Clifford algebra setting for $n=2.$ The formula (\ref{addyyeq3}), however, cannot be applied to the case due to its restriction $n>2.$ In the Clifford analysis setting, when $n>2$ the non-commutativity nature occurs, and, as a consequence, $\Gamma^*_{\xi}+\Gamma_{\xi}\ne 0.$

\end{remark}

Let $x\in {\R}^n$, $g(x): {\R}^{n}\rightarrow {\mathcal {CL}}_{0, {n-1}}$ be a Clifford-valued function.

\begin{definition}
$g(x)$ is said to be left-monogenic in domain $\Omega\subseteq {\R}^{n}$, if $g\in C^1(\Omega)$ and $Dg=0$. If $g$ is left-monogenic, then we call $\bar{D}g$ the left-derivative of $g$.
\end{definition}

Clearly, if $g$ is left-monogenic, then $g$ and all its components are harmonic.\\

Using the radial form of $D$ and formula (\ref{addyyeq3}), we know that if $P_k(r\xi)$ is a $k-$homogeneous left-monogenic function in ${\R}^n$, then $\Gamma_{\xi}P_k(\xi)=-kP_k(\xi)$ and $\Gamma^*_{\xi}P_k(\xi)=(n+k-1)P_k(\xi)$.

\begin{lemma}$^{\mbox{\scriptsize \cite{Q2}}}$\label{lemYY14}

Let $U$ be a real-valued function defined on a star-shaped open set $\Omega\subseteq {\R^n}$. If $U$ is harmonic, then there uniquely exists a left-monogenic function $F$ defined on $\Omega$ such that $Re F=U$ and $F(0)=U(0),$ given by
$$F(x)=U(x)+2 {\rm{NSc}}\int_0^1 s^{n-2}\bar{D}U(sx)x ds.$$
\end{lemma}

Using Lemma \ref{lemYY14}, the conjugate Poisson kernel in $S^{n-1}$ was first given in \cite{Brackx}. It was separately deduced in \cite{QY} based on the concept of conjugate harmonic function of the Cauchy type.  That is:

\begin{lemma}
For $\xi, \eta\in S^{n-1}$ and $r<1$. The conjugate Poisson kernel in $S^{n-1}$ is given by
\begin{equation}\label{eq25}
Q_r(\eta, \xi)=\frac{1}{\omega_{n-1}}\left[\frac{2}{|\eta-r\xi|^{n}}-\frac{n-2}{r^{n-1}}\int_0^r\frac{\rho^{n-2}}{|\eta-\rho\xi|^n}d\rho\right] r(\bar{\eta}\xi-\xi\cdot\eta)\nonumber.
\end{equation}
\end{lemma}
\bigskip

Accordingly, the Schwarz kernel in $B_n$ is defined as

\begin{eqnarray}
S_r(\eta, \xi)&=&P_r(\eta, \xi)+Q_r(\eta, \xi)\nonumber\\
&=&\frac{1}{\omega_{n-1}}\left[\frac{1-r^2}{|\eta-r\xi|^n}+(\frac{2}{|\eta-r\xi|^{n}}-\frac{n-2}{r^{n-1}}\int_0^r\frac{\rho^{n-2}}{|\eta-\rho\xi|^n}d\rho\mbox{ })r(\bar{\eta}\xi-\xi\cdot\eta)\right].\nonumber\\\nonumber
\end{eqnarray}

As in the one complex variable case, using the Schwarz kernel we can express a left-monogenic function $F$ in $B_n$ using the boundary value $u$ of the real part $U$ of $F$. That is:
\begin{equation}\label{eq26}
F(x)=F(r, \xi)=\int_{S^{n-1}}u(\eta)S_r(\eta, \xi)dS_{\eta}.
\end{equation}

\bigskip
\begin{remark}
The conjugate Poisson kernel and hence the Schwarz kernel are mixed with $2$-forms. This requires different treatment in the general Clifford algebra case in comparison with the quaternionic one. In the latter only zero- and one-forms are involved.
\end{remark}
\bigskip

We are to concern generalizations to Euclidean spaces of Fueter's Theorem on inducing quaternionic regular functions (\cite{Sce}). Obeying the same philosophy as for the one complex variable any scalar-valued function $u$ in $L^2(S^{n-1}) (n>2)$ may be expanded into a series of positive and negative spherical homogeneous monogenics (cf. page 391 \cite{Q2}) converging in the $L^2$-sense:
$$u(\xi)=\sum_{k=-\infty}^{+\infty}\frac{1}{\omega_{n-1}}
\int_{S^{n-1}}P^{(k)}(\eta^{-1}\xi)u(\eta)dS_{\eta},$$
where for the negative $k$'s the functions $P^{(k)}$ are the images of the Fueter-Sce-Qian mapping of the analytic functions $z^{k}$ in one complex variable, and for the positive $k$'s the functions $P^{(k-1)}$ are defined to be the Kelvin inversions of the corresponding $P^{(-k)}$ (\cite{Q1}). There hold $P^{(k)}(\eta^{-1}\xi)=C_{n, k}^+(\xi, \eta), k=0, 1, \cdots$ and $P^{(k)}(\eta^-\xi)=C_{n, |k|-1}^-(\xi, \eta), k=-1, -2, \cdots,$ where
\[C_{n, k}^+(\xi, \eta)=\frac{n+k-2}{n-2}C_k^{\frac{n-2}{2}}(\xi\cdot\eta)+
C_{k-1}^{\frac{n}{2}}(\xi\cdot\eta)(\bar{\eta}\xi-\xi\cdot\eta)\] and \[C_{n, |k|-1}^-(\xi, \eta)=\frac{|k|}{n-2}C_{|k|}^{\frac{n-2}{2}}(\xi\cdot\eta)-
C_{|k|-1}^{\frac{n}{2}}(\xi\cdot\eta)(\bar{\eta}\xi-\xi\cdot\eta).\]
 In particular, the Gegenbauer polynomials
$C_k^{\frac{n-2}{2}}$  satisfy $$\frac{1}{\omega_{n-1}}\frac{n+2k-2}{n-2}C_k^{\frac{n-2}{2}}=c_{n, k}P_k^n.$$

In \cite{QY} we obtain the Abelian summation representations of the Poisson and the conjugate Poisson kernels in $B_n$.

\begin{lemma}$^{\mbox{\scriptsize \cite{QY}}}$\label{lem20}

The Abel sum expansions of the Poisson and the conjugate Poisson kernels in the unit ball $B_n$ are, respectively,
$$P_r(\xi, \eta)=\frac{1}{\omega_{n-1}}\sum_{k=-\infty}^{+\infty}r^{|k|}P^{(k)}(\eta^{-1}\xi),$$
and
$$Q_r(\xi, \eta)=\frac{1}{\omega_{n-1}}
\left[\sum_{k=1}^{+\infty}\frac{k}{n+k-2}r^{k}P^{(k)}(\eta^{-1}\xi)-
\sum_{k=-\infty}^{-1}r^{|k|}P^{(k)}(\eta^{-1}\xi)\right].$$
\end{lemma}
Using Lemma \ref{lem20} and formula (\ref{eq26}), for $r<1$, we get cancelation for the negative index summation part and have
\begin{eqnarray}\label{addyeq2}
F(x)&=&\frac{1}{\omega_{n-1}}
\sum_{k=0}^{+\infty}\frac{n+2k-2}{n+k-2}r^{k}\int_{S^{n-1}}
P^{(k)}(\eta^{-1}\xi)u(\eta)dS_{\eta}\nonumber\\
&=&\frac{1}{\omega_{n-1}}\sum_{k=0}^{+\infty}\frac{n+2k-2}{n+k-2}r^{k}
\int_{S^{n-1}}C_{n,k}^+(\xi, \eta)u(\eta)dS_{\eta}.
\end{eqnarray}

\begin{theorem}
Let $F(x)\in C^1(\bar{\Omega})$. If $F(x)$ is left-monogenic in $\Omega\subseteq \R^n$, then we have
\begin{equation}\label{addeqy1}
\int_{\Omega}|\bar{D}F|^2 dV=\frac{1}{2}{\rm Sc}\int_{\partial \Omega}\bar{F}d\sigma{\bar D}F.
\end{equation}
\end{theorem}
{\bf Proof: } Using Stokes' formula, we have
$$\int_{\partial \Omega}\bar{F}d\sigma{\bar D}F=\int_{\Omega} d(\bar{F}d\sigma{\bar D}F).$$
There holds
\begin{eqnarray*}
d(\bar{F}d\sigma{\bar D}F)&=&d(\bar{F}d\sigma){\bar D}F+\bar{F}d(d\sigma{\bar D}F)\\
&=&d\bar{F}\wedge d\sigma {\bar D}F+(-1)^{n-1}\bar{F}d\sigma\wedge d({\bar D}F)\\
&=&\left[2(\bar{F}D){\bar D}F+2\bar{F}(D{\bar D}F)\right]dV\\
&=&2\left[(\bar{F}D){\bar D}F+\bar{F}(\frac{1}{4}\Delta_x F)\right]dV\\
&=&2\left[(\bar{F}D){\bar D}F\right]dV.
\end{eqnarray*}
By taking scalar part on both sides, we get (\ref{addeqy1}). This completes the proof.\\

If, in particular, $F(x)$ is left-monogenic in a ball $B_n(r)$,  then $F(x)$ has the form (\ref{addyeq2}).
Using the radial decomposition (\ref{yyeqadd1}) of $\bar{D}$, we have
\begin{eqnarray}
\bar{D}F(x)&=&\frac{1}{2\omega_{n-1}}\sum_{k=1}^{+\infty}\frac{n+2k-2}{n+k-2}r^{k-1}[ k\int_{S^{n-1}}C_{n, k}^+(\xi, \eta)u(\eta)dS_{\eta}\nonumber\\
&+&\int_{S^{n-1}}\Gamma^*_{\xi}C_{n, k}^+(\xi, \eta)u(\eta)dS_{\eta}]\bar{\xi}\nonumber\\
&=&\frac{1}{\omega_{n-1}}\sum_{k=1}^{+\infty}\frac{n+2k-2}{n+k-2}(k+\frac{n-1}{2})r^{k-1} \int_{S^{n-1}}C_{n, k}^+(\xi, \eta)u(\eta)dS_{\eta}\bar{\xi}.\nonumber
\end{eqnarray}

Using the orthogonality between $C_{n, k}^+(\xi, \eta)$and $C_{n, j}^+(\xi, \eta)$ when $k\not=j,$  we have
\begin{eqnarray}\label{addeqy20}
\int_{B_n(r)}|\bar{D}F|^2 dV&=&\frac{1}{2}{\rm Sc}\int_{S_r^{n-1}}\bar{F}d\sigma{\bar D}F\nonumber\\
&=&\frac{1}{2{\omega_{n-1}}^2}\sum_{k=1}^{+\infty}(\frac{n+2k-2}{n+k-2})^2(k+\frac{n-1}{2})
kr^{n+2k-2}\int_{S^{n-1}}\int_{S^{n-1}}u(\eta_1)
u(\eta_2)dS_{\eta_1}dS_{\eta_2}\nonumber\\
&\times &{\rm Sc}\int_{S^{n-1}}\overline{C_{n, k}^+(\xi, \eta_1)}\xi C_{n, k}^+(\xi, \eta_2)\bar{\xi}dS_{\xi}\nonumber\\
&=&\frac{1}{2}\sum_{k=1}^{+\infty}(k+\frac{n-1}{2})r^{2k+n-2}c_{n,k}\int_{S^{n-1}}\int_{S^{n-1}}u({\eta_1})u({\eta_2})P^n_k({\eta_1}\cdot{\eta_2})dS_{{\eta_1}}dS_{{\eta_2}}\nonumber\\
&=&\frac{1}{2}\sum_{k=1}^{+\infty}(k+\frac{n-1}{2})r^{2k+n-2}\sum_{j=1}^{a^n_k}|b_j|^2.
\end{eqnarray}

Taking limit $r\rightarrow 1^-$ on both sides of (\ref{addeqy20}), we have
$$\int_{B_n} |\bar{D}F|^2 dV=\frac{1}{2}{\rm Sc}\int_{S^{n-1}}\bar{F}d\sigma{\bar D}F=\frac{1}{2}\sum_{k=1}^{+\infty}(k+\frac{n-1}{2})\sum_{j=1}^{a^n_k}|b_j|^2\not=\int_{B_n}|\bigtriangledown U|^2dV.$$
Then formula (\ref{eq23}) is proved.

In summary, in the Clifford algebra setting of the n-dimensional Euclidean space we have

\begin{theorem}
Let $F$ be left-monogenic in $B_n$ $(n>2)$ and $U$ its scalar part with square-integrable boundary value $u$ in $S^{n-1}.$ Then there holds
\begin{eqnarray*}
A(u)&:=&\frac{1}{\omega_{n-1}}\int_{S^{n-1}}\int_{S^{n-1}}\frac{|u({\eta_1})-u({\eta_2})|^2}{|{\eta_1}-{\eta_2}|^n}
 dS_{{\eta_1}}dS_{{\eta_2}}\\
&=&\int_{B_n}|\bigtriangledown U|^2dV\\
&=&\sum_{k=1}^{+\infty}k\sum_{j=1}^{a^n_k}|b_j|^2
\end{eqnarray*}

and
\begin{eqnarray*}
&&\int_{B_n}|\bar{D}F|^2 dV\\
&=&\frac{1}{2}{\rm Sc}\int_{S^{n-1}}\bar{F}d\sigma{\bar D}F\\
&=&\frac{1}{2}\sum_{k=1}^{+\infty}(k+\frac{n-1}{2})\sum_{j=1}^{a^n_k}|b_j|^2,
\end{eqnarray*}
where $b_j$'s are the coefficients of the Fourier-Laplace series of $u$ which are given by formula (\ref{eqY14}).

\end{theorem}

\begin{remark}
We note that in the present section, \S 4, the $n=3$ case corresponds to ${\mathcal {CL}}_{0, {n-1}}={\mathcal {CL}}_{0, 2}=\H,$ the quaternions. The theory in \S 4, however, does not give the one established in \S 3. In fact, the domain of $u$ in \S 4 is $S^{2}$ contained in ${\mathbb R}^3,$ while the domain of $u$ in \S 3 is $S^{3}$ contained in ${\mathbb R}^4.$ The quaternionic case can take the Cauchy-Riemann equation approach which is different and simpler.
\end{remark}

{\bf Acknowledgement}. We wish to sincerely thank Huaying Wei who provided reference \cite{Ahlfors} with benificial instructions.

\end{document}